%% file: ms.tex
\documentclass[11pt]{article}
\usepackage[utf8]{inputenc}
\usepackage[T1]{fontenc}
\usepackage{lmodern}

\usepackage{graphicx}
\usepackage{color}
\usepackage{hyperref}

\usepackage{amsmath,amsfonts,amssymb,amsthm}
\theoremstyle{plain}
\newtheorem{thm}{Theorem}
\newtheorem{lemma}{Lemma}

\title{Nested distance for stagewise-independent processes}
\author{Filipe Goulart Cabral, Bernardo Freitas Paulo da Costa}
\date{\today}

\begin{document}

\maketitle

\begin{abstract}
We prove that, for two discrete-time stagewise-independent processes with a stagewise metric,
the nested distance is equal to the sum of the Wasserstein distances
between the marginal distributions of each stage.
\end{abstract}

\tableofcontents

\section{Introduction}

% Problema: discretizar os cenários (ou reduzir)
%   Métodos típicos: Sampling ; Wasserstein e ND
%   Caso em questão: SDDP <=> stagewise indep
%     Resultado: neste caso, ND = Soma das W

An usual approach when solving a multi-stage stochastic programming problem is approximating the underlying probability distribution by a \emph{scenario tree}.
Once we obtain this approximation, the resulting problem becomes a deterministic optimization problem with an extremely large number of variables and constraints.
Most of the standard algorithms for solving it depend on the convexity properties of the objective function and constraints:
for example the Cutting Plane (or L-Shaped), Trust Region, and Bundle methods, for the two-stage case;
and Nested Cutting Plane, Progressive Hedging and SDDP, for the multi-stage case,
see~\cite{birge2011introduction,ski2003stochastic}.
% Refs

Although all those methods are very effective,
in practice their performance depend heavily on the size of the scenario tree.
Ideally, we seek the best approximation with the least number of scenarios,
since a large number of scenarios can really impact computational time.
There are two main available techniques for the scenario generation:
those based on sampling methods (like Monte Carlo) and those based on optimal scenario generation using probability metrics.
The later is the subject of this article.

In particular, the widely used SDDP method~\cite{pereira1991multi} depends on an additional hypothesis about the underlying uncertainty:
the \emph{stagewise independence property}~\cite{shapiro2011analysis}.
We postpone a precise definition to section~\ref{sec:swi},
% but we note that the success of the SDDP algorithm comes from a significant reduction of the computational cost
% that is obtained thanks to this hypotesis.
but we note that such property is crucial for the significant  computational cost reduction of the SDDP algorithm.
Our aim here is to show that, analogously, the stagewise independence
allows for a similar reduction of the computation time in a related problem regarding the optimal scenario generation technique.

% Re-escrever para começar direto no success of SDDP??
% Emphasize the importance of the stagewise independence in this article.

\subsection{Optimal scenario generation}
%In this section, , since it is the main technique of this article.
% and how the  paper deals with the optimal scenario generation for the multi-stage and stagewise independent setting.
%The main result relating the

Optimal scenario generation aims to approximate,
in a reasonable sense for stochastic optimization,
a general probability distribution by a discrete one with a fixed number of support points.
There are several probability metrics that could be used as objective functions for the optimal probability discretization:
an extensive list with 49 examples can be found in~\cite[\textsection 14]{DezaDeza201410}.
Among all those probability metrics, the Wasserstein distance stands out as a convenient one, \cite{pflug2001scenario,heitsch2003scenario,dupavcova2003scenario},
since under some mild regularity conditions it provides an upper bound
on the absolute difference between the optimal values
of a two-stage stochastic optimization problem with distinct probability distributions~\cite[section 2 -- page 242]{romisch1991stability}.

A generalization of the Wasserstein distance for the multi-stage case
which has an analogous bound for the difference between optimal values is the Nested Distance
developed in~\cite{pflug2009version,pflug2012distance,pflug2014multistage}.
% Ref
The standard algorithm for evaluating it is based on a dynamic programming problem~\cite{KovacevicP15},
whose intermediate subproblems are similar to (conditional) Wasserstein distances.
However, this is prohibitive for general distributions,
since the number of intermediate subproblems grows with the number of scenarios,
which is generally exponential on the number of stages.
For practical applications, optimal scenario generation using the Nested Distance is therefore very limited.

In this paper we show that, for stagewise independent distributions,
it is possible to reduce dramatically the computational burden of evaluating the Nested Distance.
%The idea is exploring symmetries of the corresponding probability measures in the dynamic programming formulation.
Actually, we obtain a stronger result relating the Nested and Wasserstein distances:
We prove that the Nested distance is equal to the sum of the Wasserstein distances
between the marginal distributions of each stage.
In particular, the number of subproblems required for evaluating the Nested Distance is now equal to the number of stages,
and each subproblem can be solved independently and very effectively by calculating the Wasserstein distances.
This result supports a new scenario reduction method that preserves the stagewise independence property,
which was compared to the standard Monte Carlo approach in~\cite[\textsection Appendix A]{shapiro2015guidelines}.

\bigskip

\subsection{Organization}
In this paper, we focus on the case of evaluating the nested distance between discrete-time, discrete-state stochastic processes.
This is not very restrictive, since in most cases the best we can do is producing very large samples and computing the nested distance from them, due to the complexity of the nested distance formula.

In section~\ref{sec:w_nd}, we will review the definitions of the Wasserstein distance and the Nested distance.
We also present the usual tree representation of discrete-time stochastic process
as a motivation for a matrix representation of the linear problems defining both distances.

Then, in the following section, we recall the definition of stagewise independence for processes,
and observe how this assumption simplifies the trees corresponding to them.
This suggests an analog simplification for the Nested Distance picture,
which we will prove correct in section~\ref{sec:3stage} in the fundamental 3-stage setting.

Finally, in section~\ref{sec:multistage},
we recall the different equivalent linear programming formulations of the Nested distance.
Then, we define the subtree distance, which will be,
along with the intuition developed in the 3-stage case,
the fundamental tool to proving our result.

\bigskip

%Agradecimentos
%Joari, Tito, Erlon,  ONS.

We thank professor Tito Homem-de-Mello, Universidad Adolfo Ibanez, and Erlon C. Finard, Federal University of Santa Catarina, for the enlightening discussion occurred on the XIV International Conference in Stochastic Programming, which have encouraged us for writing this paper.
We would like to show our gratitude to Joari P. da Costa, Brazilian Electrical System Operator (ONS), for the assistance and comments that greatly improved the manuscript.
We are also grateful to Alberto S. Kligerman, ONS, for the opportunity to conduct this research.

\input wasserstein_nd

\input swi

\input 3stage

\input multistage

\input example

\bibliographystyle{abbrv}
\bibliography{refs.bib}

\end{document}

%% file: wasserstein_nd.tex
\section{Wasserstein and Nested distances}
\label{sec:w_nd}

Before presenting the Nested distance,
we review the definition of Wasserstein distance and some of its properties, following the notation of ~\cite{KovacevicP15}.
This motivates the introduction of the Nested distance
and it will also be used in the conclusion.

\subsection{Wasserstein distance}
We start with a very general definition.
Let $(\Xi,\mathcal{F}, P)$ and $(\Xi,\mathcal{G}, Q)$ be two probability spaces and $d: \Xi \times \Xi \rightarrow \mathbb{R}$ be a distance function.
The \emph{Wasserstein distance} of order $p \geq 1$ between both probability spaces, denoted by $d_{W,r}(P,Q)$,
is the optimal value of the optimization problem
\begin{equation} \label{eq:wasserstein}
\begin{array}{rll}
\min_{\pi}    & \displaystyle \int d^p(\xi,\zeta) \cdot \pi(d\xi,d\zeta) \\
\textrm{s.t.} & \pi(M \times \Xi) = P(M)
              & \text{for all } M \in \mathcal{F}, \\
              & \pi(\Xi \times N) = Q(N)
              & \text{for all } N \in \mathcal{F},  \\
              %& \text{$\pi$ a probability on $(X \times X, \mathcal{F} \otimes \mathcal{F})$.}
\end{array}
\end{equation}
where the minimum on \eqref{eq:wasserstein} is among all probability measures $\pi$
on the product space $(\Xi \times \Xi, \mathcal{F} \otimes \mathcal{G})$.

If $P = \sum_{i = 1}^m p_i \delta_{\xi^i}$ and $Q = \sum_{j = 1}^n q_j \delta_{\zeta^j}$
are two discrete probability measures,
then the Wasserstein distance can be computed by the following linear program:
\begin{equation}\label{eq:wassersteins_lp}
\begin{array}{rll}
\displaystyle \min_{\pi \in \mathbb{R}^{m \times n}}
              & \displaystyle \sum_{i,j} d_{i,j}^r \cdot \pi_{i,j}  \\
\textrm{s.t.} & \sum_{j} \pi_{i,j} = p_i
              & \text{for all $i \in \{1,\dots, m\}$}, \\
              & \sum_{i} \pi_{i,j} = q_j
              & \text{for all $j \in \{1,\dots, n\}$},  \\
              & \sum_{i,j} \pi_{i,j} = 1, \\
              & \pi_{i,j} \geq 0,
              & %(i = 1,\dots, m, \ j = 1,\dots,\nu),
\end{array}
\end{equation}
where $p_i, q_j$ are the corresponding probabilities of the outcomes $\xi^i, \zeta_j$,
and the linear coefficients from the objective function~$d_{i,j}^r := d^r(\xi^i, \zeta^j)$ are the distances between those outcomes.
Note that the constraint~$\sum_{i,j} \pi_{i,j} = 1$ is redundant
since it follows from any of the first two sets of constraints in~\eqref{eq:wassersteins_lp},
\[
\sum_{i = 1}^m \sum_{j = 1}^n \pi_{i,j} = \sum_{i = 1}^m p_i = 1.
\]

\begin{figure}[ht]
  \centering
  \includegraphics[width=0.4\textwidth]{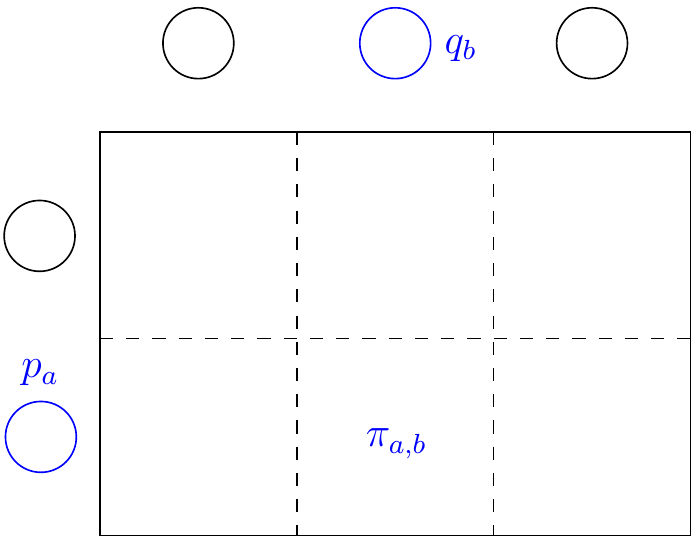}
  \caption{A graphical representation of the tableau for $\pi_{a,b}$ calculating the Wasserstein distance,
  and its relation to the marginal probabilities $p$ and $q$.}
  \label{fig:wd}
\end{figure}

The constraints from the Wasserstein distance problem~\eqref{eq:wasserstein} and~\eqref{eq:wassersteins_lp} impose that the joint probability distribution $\pi$ on the product space $(\Xi \times \Xi, \mathcal{F} \otimes \mathcal{G})$ must have $P(\cdot)$ and $Q(\cdot)$ as marginal distributions.
Moreover, the objective function guarantees that the optimal joint probability $\pi^*$ induces the least expected value $\mathbb{E}_{\pi^*}[d^r(\xi,\zeta)]$ for the distance function $d^r$.
Figure~\ref{fig:wd} illustrates the constraints from~\eqref{eq:wassersteins_lp} in terms of a joint probability table.
%An illustration of the constraints from~\eqref{eq:wassersteins_lp} is given in figure~\ref{fig:wasserstein} in terms of a joint probability table.
Each row label $p_i$ and column label $q_j$ correspond to the probability of a given outcome $i$ and $j$% from $P$ and $Q$
, respectively,
and the cells represent the values $\pi_{i,j}$ of the associated joint probabilities.
So, the column sum $\sum_{j} \pi_{i,j}$ must be equal to the row probability $p_i$ as well as the row sum $\sum_{i} \pi_{i,j}$ must be equal to the column probability $q_j$.
The generalization of such joint probability table is useful to visualize the constraints from the nested distance case in section~\ref{sec:nested_distance}.

Another instructive interpretation of problems~\eqref{eq:wasserstein} and~\eqref{eq:wassersteins_lp} is in terms of optimal transportation.
The constrains of \eqref{eq:wassersteins_lp} can be seen as the transportation of goods from $m$ sources to $n$ destinations.
The sources are indexed by $i$ and each of them have $p_i$ goods available (we can imagine $p_i$ as a fraction of a total number).
The destinations are indexed by $j$ and each of them have $q_j$ of demand.
So, the decision variable $\pi_{i,j}$ corresponds to the proportion of goods to send from source $i$ to destination $j$ and the parameter $d_{i,j}^r$ is the unit cost of such shipment.
Therefore, problem~\eqref{eq:wassersteins_lp} is an optimal transportation problem whose objective function is minimizing the overall cost.
Problem~\eqref{eq:wasserstein} follows the same reasoning, but for the general case which includes continuous mass transportation from a region, e.g., subset of $\mathbb{R}^n$, to another.

\subsection{Notation for probability trees}

Discrete-time, discrete-state stochastic processes have a one-to-one correspondence with probability trees.
The latter is much more intuitive and easy to deal with, so we describe the underlying uncertainty by probability trees.
Since we're always dealing with two parallel objects,
we'll introduce a notation that hopefully simplifies our discussion
and keeps this parallel as clear as possible.

\begin{figure}[ht]
  \centering
  \includegraphics[width=0.9\textwidth]{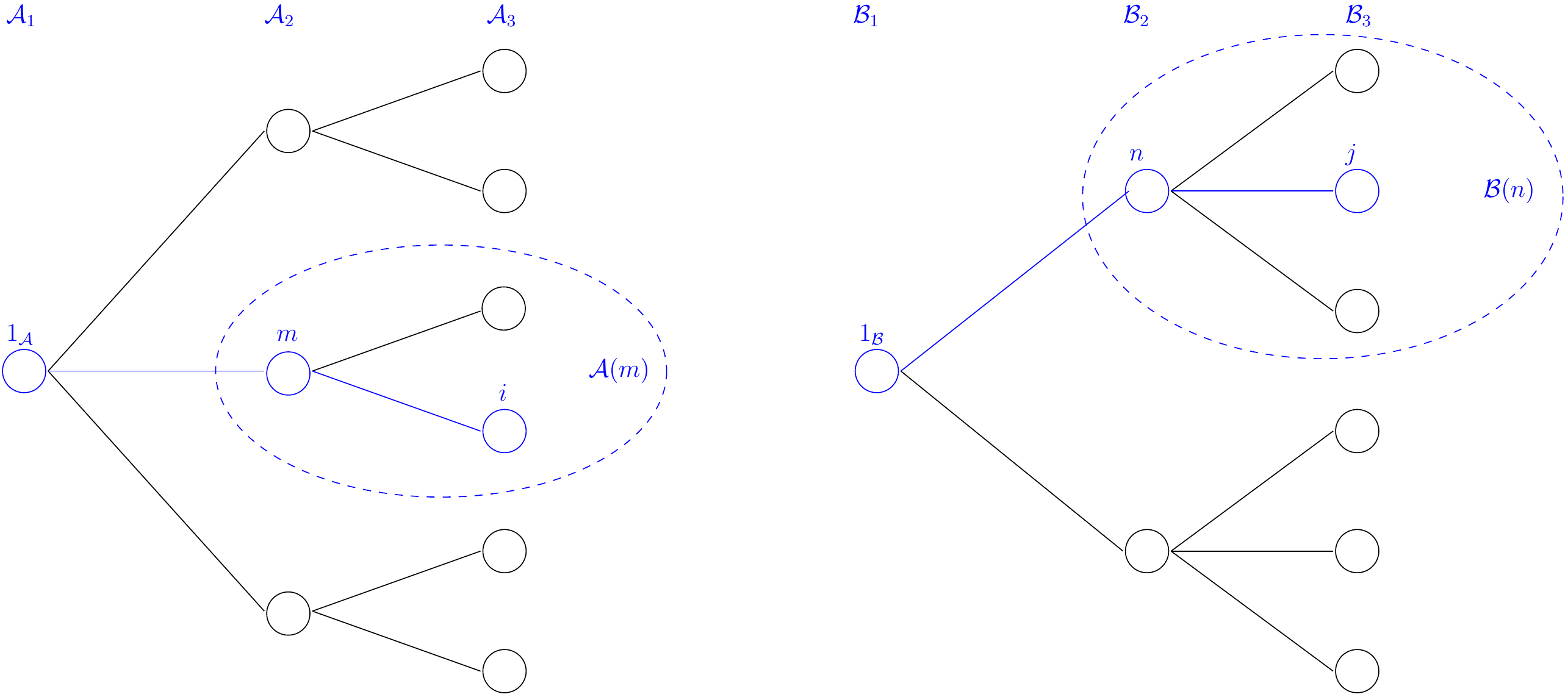}
  \caption{The trees $\mathcal{A}$ and $\mathcal{B}$ in a 3-stage setting}
  \label{fig:3stagetree}
\end{figure}

Our probability trees are called $\mathcal{A}$ and $\mathcal{B}$,
and both have $T$ stages, see figure~\ref{fig:3stagetree} for a three stage example.
Those probability trees have subsets $\mathcal{A}_t$ and $\mathcal{B}_t$
which correspond to nodes at stage $t$.
A general \emph{leaf} will be denoted respectively by $i$ and $j$,
a general \emph{node}, $k$ and $l$,
and (when needed) another general node/leaf, $m$ and $n$.
The roots are respectively denoted by $1_\mathcal{A}$ and $1_\mathcal{B}$,
but when there's no ambiguity we'll just refer to them as $1$.
The subtree rooted at node $k$ will be $\mathcal{A}(k)$
and the notation $i \succ k$ asserts that the \emph{leaf} $i$ belongs to $\mathcal{A}(k)$.
Similarly, we have $\mathcal{B}(l)$ and $j \succ l$.
Note that a leaf node $i \in \mathcal{A}_T$ defines uniquely a path from the root node $1_{\mathcal{A}}$ to the leaf~$i$.
Such path is called a scenario and is illustrated in blue in figure~\ref{fig:3stagetree}.
The node~$m$ is called the \emph{successor} node of~$k$, denoted by~$m \in k_+$, if~$k$ is the immediate predecessor of~$m$ in the path that connects the root~$1_{\mathcal{A}}$ to $m$.
The same concept also follows for the tree~$\mathcal{B}$.

%The probabilities of $\mathcal{A}$ and $\mathcal{B}$ are denoted by $P(\cdot)$ and $Q(\cdot)$, respectively.
%We shall describe only the probabilities of $\mathcal{A}$, since the other case is analogous.

The probability of~$\mathcal{A}$ is described by the probability function~$P(\cdot)$ defined at the leaf nodes (or scenarios) of~$\mathcal{A}$ and extended to a general node~$k \in \mathcal{A}_t$ by summing up the probabilities of each leaf node~$i \in \mathcal{A}_T$ descendant from~$k$:
\begin{equation}\label{eq:node_prob}
  P(k) = \sum_{i \succ k} P(i).
\end{equation}
The conditional probability~$P(i \mid k)$ of a leaf~$i$ conditioned to a given node~$k$ is therefore
\[
P(i \mid k) = %\frac{P(i)}{P(k)}.
\left\{
\begin{array}{ll}
\frac{P(i)}{P(k)}
%P(i)/P(k)
    & \text{, if $i \succ k$} \\
0   & \text{, otherwise}
\end{array}
 \right..
\]
Analogously to~\eqref{eq:node_prob}, conditional probability is extended to a general node~$m$ by the leaf nodes:
\begin{equation}\label{eq:node_condprob}
P(m \mid k) = \sum_{i \succ m} P(i \mid k).
\end{equation}
A relation between the probability of a given node $k$ and its successors, $m \in k_+$, is given by the following identity:
\begin{equation}\label{eq:node_identity}
P(k) = \sum_{m \in k_+} \sum_{i \succ m} P(i) = \sum_{m \in k_+} P(m).
\end{equation}
Formula \eqref{eq:node_identity} will be useful in the proof of the main result of this paper and is analogous to conditional probabilities.
The same concepts also apply to the tree~$\mathcal{B}$ and probability function~$Q(\cdot)$.

\subsection{The nested distance} \label{sec:nested_distance}

Let $\mathcal{A}$ and $\mathcal{B}$ be two probability trees with the same number of stages as above.
The nested distance of order $p \geq 1$ between $\mathcal{A}$ and $\mathcal{B}$ is defined as the optimal value of the optimization problem
\begin{equation} \label{eq:ND}
\begin{array}{rll}
\displaystyle \min_{\pi} & \displaystyle \sum_{i,j} d_{i,j}^p  \cdot \pi_{i,j} \\[1ex]
\textrm{s.t.} & \text{for all $k \in \mathcal{A}_t$ and $l \in \mathcal{B}_t$:} \\
              & \sum_{j \succ l} \pi(i, j \mid k, l) = P(i \mid k)
              & \text{for all $i \succ k$}, \\
              & \sum_{i \succ k} \pi(i, j \mid k, l) = Q(j \mid l)
              & \text{for all $j \succ l$}, \\[1ex]
              & \sum_{i,j} \pi_{i,j} = 1, \\
              & \pi_{i,j} \geq 0,
              & \text{for all leaves $i \in \mathcal{A}_T$, $j \in \mathcal{B}_T$},
\end{array}
\end{equation}
where the minimum on \eqref{eq:ND} is among all discrete probability measures $\pi$ on the product tree $(\mathcal{A} \times \mathcal{B})$.
Each part of problem \eqref{eq:ND} is described in details below:
\begin{itemize}
  \item the equality constraint
    \begin{equation} \label{eq:NDcst1}
      \sum_{j \succ l} \pi(i, j \mid k, l) = P(i \mid k)
    \end{equation}
    enforces that the marginal distribution of $\pi(i,j \mid k, l)$ on the probability tree $\mathcal{B}$ must be equal to the conditional probability $P(i \mid k)$ from $\mathcal{A}$.
    % in the left-hand side is the sum of conditional probabilities $\pi(i,j \mid k, l)$ among all leaves $j \in \mathcal{A}_T$ which have the node $l \in \mathcal{A}_t$ as its predecessor.
    % This sum represents the marginal distribution of $\pi(i,j \mid k, l)$ on the probability tree $\mathcal{B}$ and the constraint \eqref{eq:NDcst1} enforce that it must be equal to the conditional probability $P(i \mid k)$ from $\mathcal{A}$.
    Note that problem~\eqref{eq:ND} considers all possible combinations of the marginal constraint~\eqref{eq:NDcst1}, i.e., it considers one constraint of type~\eqref{eq:NDcst1} for each  stage~$t$,
    each pair of nodes~$(k,l) \in \mathcal{A}_t\times \mathcal{B}_t$
    and each leaf node~$i \in \mathcal{A}_T$ descendent from~$k$.
    The same comments apply to the marginal constraint with the conditional probability $Q(j \mid l)$ in right-hand side;

  \item The objective function is given by the sum $\sum_{i,j} d_{i,j}^p \cdot \pi_{i,j}$ among all leaf nodes $i \in \mathcal{A}_T$ and $j \in \mathcal{B}_T$.
    This sum represents the expected value of a distance function $d_{i,j}^p := d^p(\xi^i,\zeta^j)$ between observations along pairs of scenarios.
    The main example of distance $d^p(\xi^i,\zeta^j)$ is the (weighted) \textit{stagewise distance}:
    \begin{equation}\label{eq:scen_distance}
      d^p(\xi^i, \zeta^j)
        = \sum_{t = 1}^T w_t \cdot d_t^p (\xi_t^i, \zeta_t^j).
    \end{equation}
    The scenario distance \eqref{eq:scen_distance} is crucial for the main result of this paper.

  \item The non-negativity constraints and sum to one constraint,
    \[
      \pi_{i,j} \geq 0, \quad  \sum_{i,j} \pi_{i,j} = 1,
    \]
    ensure that each feasible solution from \eqref{eq:ND} is a probability distribution on the tree product $\mathcal{A}\times\mathcal{B}$.
    We emphasize that constraints relating total and conditional probabilities for $\pi$ on $\mathcal{A} \times \mathcal{B}$ are implicit, that is,
    \begin{itemize}
        \item $\pi(i,j \mid k,l) = \pi_{i,j}/\pi_{k,l}$;
        \item $\pi_{k,l} = \sum_{(i,j)\succ (k,l)} \pi_{i,j}$.
    \end{itemize}
\end{itemize}

\begin{figure}[ht]
	\centering
	\includegraphics[width=0.9\textwidth]{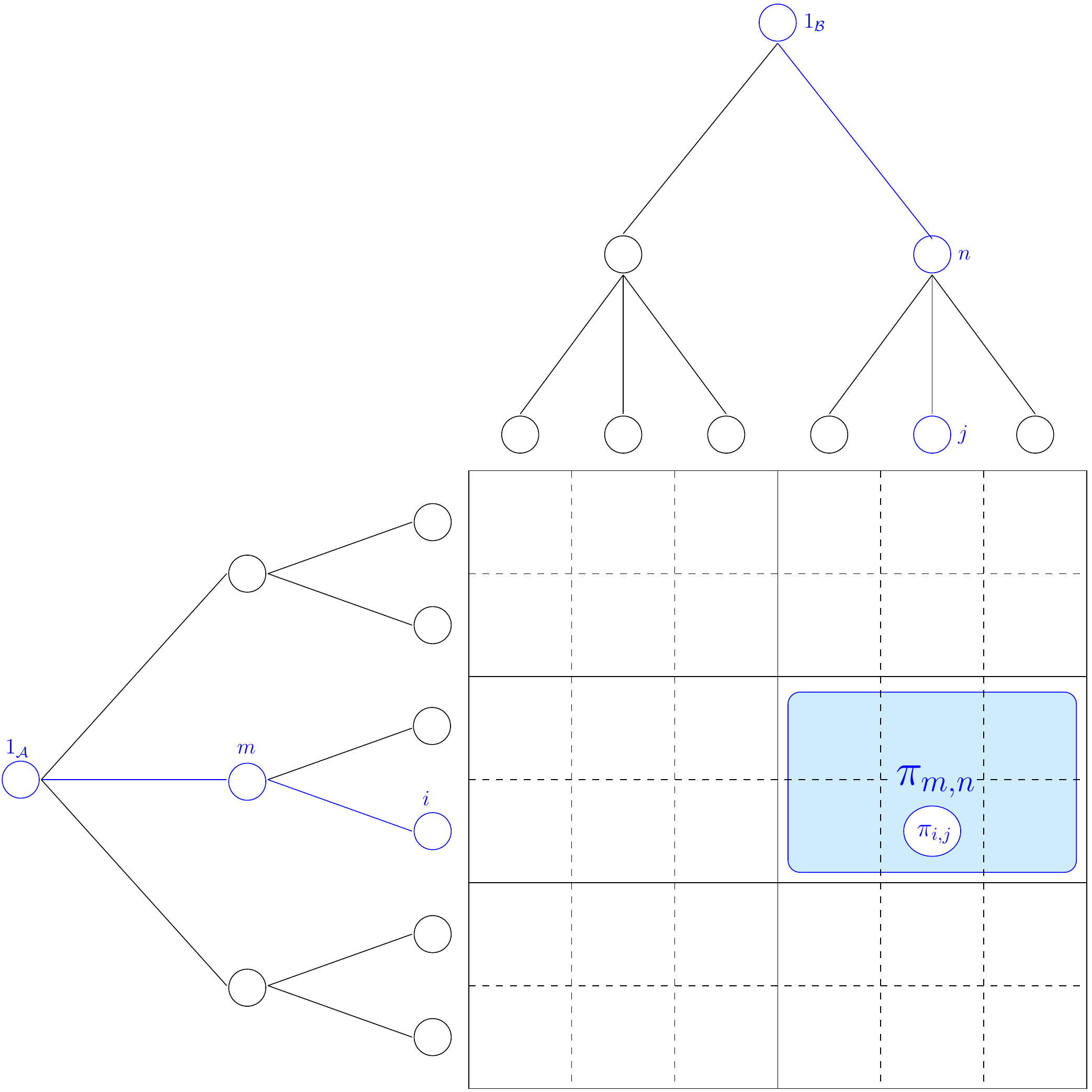}
	\caption{The trees $\mathcal{A}$ and $\mathcal{B}$ in a 3-stage setting}
	\label{fig:3stagepi}
\end{figure}

We illustrate in figure~\ref{fig:3stagepi},
inspired from~\cite[Figure 2.8]{pflug2014multistage},
the constraints of a nested distance problem using a three stage example.
The probability trees $\mathcal{A}$ and $\mathcal{B}$ are represented on the left and top side of the main block, respectively.
Each cell delineated by dashed lines corresponds to a joint probability~$\pi_{i,j}$ of leaf nodes~$i \in \mathcal{A}_3$ and~$j \in \mathcal{B}_3$, and the sum of probabilities within a solid line block corresponds to a probability~$\pi_{m,n}$ of second stage nodes $m \in \mathcal{A}_2$ and $n \in \mathcal{B}_2$.
So, the ratio of each cell to the associated block is the conditional probability~$\pi(i,j \mid m,n)$.
The sum of ratios along a line of a block, which correspond to sum along the leaves $j$ descendants from $n$, leads to the left-hand side of the marginal constraint with $P(i \mid m)$ in the right-hand side.
The marginal constraint with $Q(j \mid n)$ is obtained if we sum those ratios along a column.

% TODO
% Comment about other formulations that will be useful in the general proof.

% vim:set spelllang=en:

%% file: swi.tex
\section{Stagewise independent processes and trees}
\label{sec:swi}

A discrete-time stochastic process $(\xi_t)_{t=1}^T$ is \emph{stagewise independent} (SWI)
if the random variable $\xi_{t+1}$ is independent of its entire past
$(\xi_1, \xi_2, \ldots, \xi_t)$.

This is a very strong restriction on the random process underlying the optimization problem,
and not necessarily a realistic asumption.
However, it is very often possible to reformulate the process $(\xi_t)_{t=1}^T$
in terms of another, stagewise independent process $(\tilde\xi_t)_{t=1}^T$,
inducing an equivalent optimization problem.
% TODO: cite
The great advantage of SWI is that the resulting optimization problem
can be solved with much more efficient algorithms.
% TODO cite sddp

% TODO: think...
% The simplest way of constructing such process is through a sequence of
% independent random variables.
% More generally, (define \otimes here??)
\subsection{SWI trees as product of trees}

\begin{figure}[ht]
  \centering
  \includegraphics[width=0.9\textwidth]{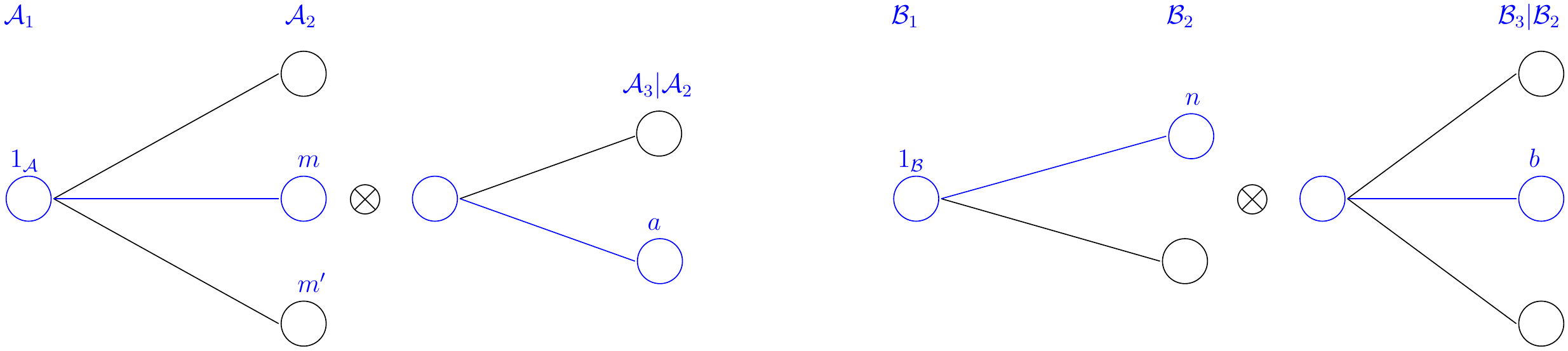}
  \caption{The SWI representation of trees $\mathcal{A}$ and $\mathcal{B}$}
  \label{fig:indeptrees}
\end{figure}

Since the probability of the events of $\xi_{t+1}$ do not depend on their past,
we can use a condensed representation of the probability tree
corresponding to $(\xi_t)_{t=1}^T$.
In the tree model, stagewise independence implies that
every node in stage $t$ has exactly the \emph{same} descendants in stage $t+1$,
with exactly the same probability.
This yields a symmetrical tree,
and indeed both trees in figure~\ref{fig:3stagetree} could correspond to SWI process,
as depicted in figure~\ref{fig:indeptrees}.

Since probability trees correspond to stochastic processes,
we can define a product operation between two trees $\mathcal{A}'$ and $\mathcal{A}''$:
The process corresponding to $\mathcal{A} = \mathcal{A}' \otimes \mathcal{A}''$
is given by realizations
\[
  a = a' \otimes a'' = (a'_1, a'_2, \ldots, a'_{T'}, a''_1, a''_2, \ldots, a''_{T''})
\]
for every pair of leaves $a'$, $a''$ of $\mathcal{A}'$ and $\mathcal{A}''$,
respectively,
with probability $P(a) = P'(a') \times P''(a'')$.
This operation establishes a bijection between the subtree
$\mathcal{A}(a')$ rooted at $a'$ (as an internal node of $\mathcal{A}$)
and the tree $\mathcal{A''}$.
In particular, it also give bijections between two subtrees of $\mathcal{A}$.

\begin{figure}[ht]
  \centering
  \includegraphics[width=0.9\textwidth]{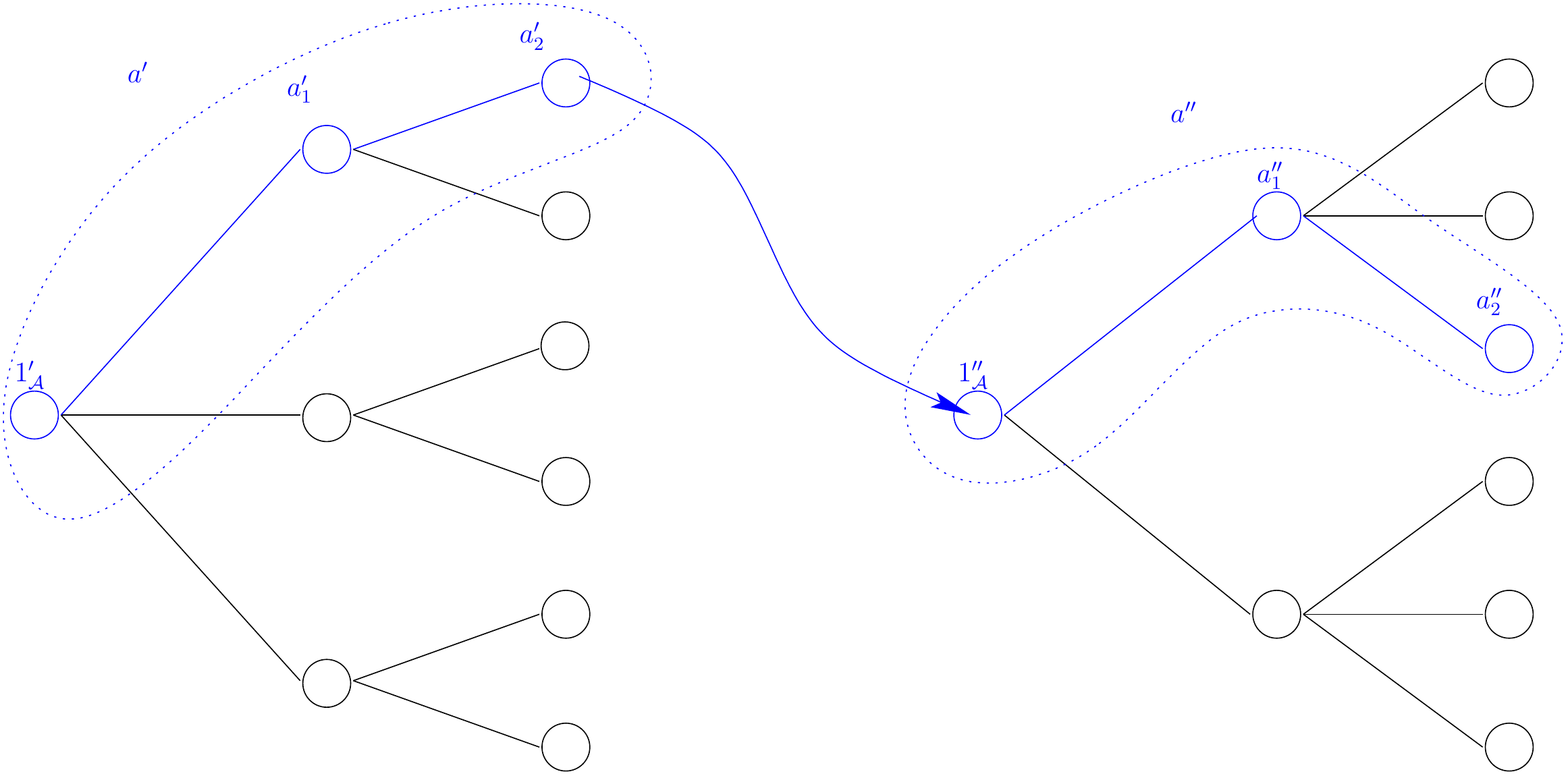}
  \caption{The product of trees $\mathcal{A'}$ and $\mathcal{A''}$}
  \label{fig:tree_prod}
\end{figure}

Visually, the $\otimes$ symbol corresponds to the operation of attaching
the probability tree on its right to \emph{every} leaf in the probability tree on the left,
with probabilities induced by multiplication as shown in figure~\ref{fig:tree_prod}.
For example, if the probability of node $a'$ is $0.4$, that of $a''$ is $0.2$,
the probability of the node $a = a' \otimes a''$ will be $0.08$.

Observe that stagewise independence is stronger than
requiring the same ``topological'' structure for every node in a given state.
Indeed, SWI also implies the equality between the ``transition probabilities''
from each node in stage $t$ to the corresponding descendant in stage $t+1$.
Therefore, when we draw a tree as a product of trees, as above,
we are asserting not only a regularity for the state-space of the process,
but also that some conditional probabilities are equal.
If these conditions were met, we could draw trees $\mathcal{A}$ and~$\mathcal{B}$
as in figure~\ref{fig:indeptrees}.

\medskip

% Motivar a decomposição da ND em SW-WD

\begin{figure}[ht]
  \centering
  \includegraphics[width=0.9\textwidth]{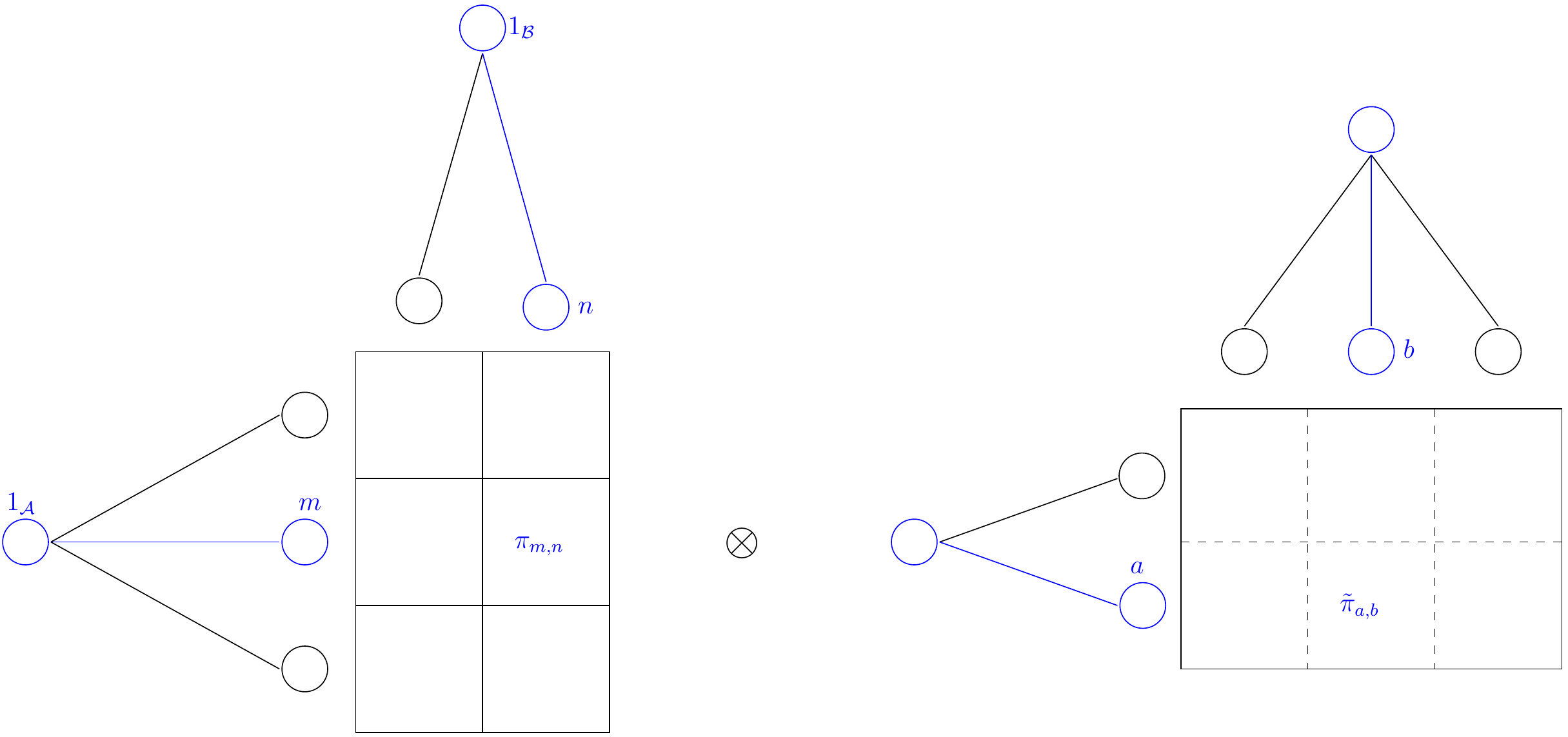}
  \caption{The tableaux for nested distance between SWI trees $\mathcal{A}$ and $\mathcal{B}$}
  \label{fig:piindep}
\end{figure}

The simple structure of SWI trees suggests that
we could be able to split the nested-distance tableau
of figure~\ref{fig:3stagepi} also as a product,
as depicted in figure~\ref{fig:piindep}.
This amounts to reducing the computation of the nested distance
as two independent subproblems of Wasserstein distances.

In the next section, we will prove that this is indeed the case for 3-stage trees,
and in the following one we will generalize, by induction,
for arbitrary stagewise independent trees.

% vim:set spelllang=en:

%% file: 3stage.tex
\section{The 3-stage case}
\label{sec:3stage}

There are two main inspirations for our result,
which are also helpful guides for the proof strategy.
The first one comes from the dynamic programming approach to multi-stage stochastic optimization problems,
where the calculations are done at the last stage $T$ and then recursively used in previous stages
until reaching the first one.
Since, as it happened there, the recursive strategy for the general $T$-stage problems already appears in the 3-stage case,
we start with this case in order to present a lighter notation.

The second one arises from the ``decomposable'' nature of many multi-stage objective functions
as \emph{sums} of stage-wise costs.
Even if the recursive formulation we obtain as an intermediate step is independent from this second hypothesis,
the final result is not, and in this section we'll exploit it early on,
in an attempt to motivate more clearly the steps we take.

\subsection{A simplified notation}

Since we don't need to deal with an arbitrary number of stages,
we will always use in this section the convention that $k$ and $l$ are nodes from the second stage,
respectively predecessors of $i$ and $j$ when needed.
For example, in $\sum_i g(k)$, $k$ is to be understood as denoting $i$'s predecessor.

\subsection{Objective function decomposition}

The distances $d_{i,j} = d(\xi^i, \zeta^j)^p$ between scenarios $\xi^i$ and $\zeta^j$ are sums of stage-wise distances:
$d_X(\xi^i_1, \zeta^j_1)^p + d_X(\xi^i_2, \zeta^j_2)^p + d_X(\xi^i_3, \zeta^j_3)^p$,
which we abbreviate as $\delta_{1,1} + \delta_{k,l} + \delta_{i,j}$.
\ (Note that if both $(i,j)$ and $(i',j')$ are successors of $(k,l)$,
their second stage restrictions are equal: $\xi^i_2 = \xi^{i'}_2$ and $\zeta^j_2 = \zeta^{j'}_2$,
so $\delta_{k,l}$ is unambiguously defined).

We rewrite the objective function in ``stage-wise terms'', and split:
\begin{align*}
  f(\pi) & = \sum_{i,j} \pi_{i,j} d_{i,j} = \sum_{i,j} \pi_{i,j} \big(\delta_{1,1} + \delta_{k,l} + \delta_{i,j}\big) \\
         & = \left(\sum_{i,j} \pi_{i,j} \delta_{1,1}\right)
                + \left(\sum_{i,j} \pi_{i,j} \delta_{k,l}\right)
                + \left(\sum_{i,j} \pi_{i,j} \delta_{i,j}\right)
\end{align*}
The first term is simply $\delta_{1,1}$, summing all the $\pi_{i,j}$.
Rearranging $\sum_{i,j}$ into $\sum_{k,l} \sum_{(i,j)\succ(k,l)}$, the third term becomes
\[
  \sum_{k,l}\sum_{(i,j)\succ(k,l)} \pi_{i,j} \delta_{i,j}.
\]
The second term is more involved, since it both needs to rearrange and sum the $\pi_{i,j}$:
\[
  \sum_{i,j} \pi_{i,j} \delta_{k,l}
  = \sum_{k,l} \sum_{(i,j)\succ(k,l)} \pi_{i,j} \delta_{k,l}
  = \sum_{k,l} \delta_{k,l} \sum_{(i,j)\succ(k,l)} \pi_{i,j}
  = \sum_{k,l} \delta_{k,l} \pi_{k,l}.
\]
where we defined $\pi_{k,l} := \sum_{(i,j)\succ(k,l)} \pi_{i,j}$,
as the total probability of the internal node pair $(k,l)$.
This is consistent with the analogous ``first-to-second-stage equation''
$\pi_{1,1} = \sum_{(k,l)\succ(1,1)} \pi_{k,l}$
which gives the total probability of the tree.

Putting it all together, we obtain two new equivalent expressions for the objective function:
\begin{align}
  \label{eq:obj_f_by_stage}
  f(\pi) & = \delta_{1,1}
                + \left(\sum_{k,l} \delta_{k,l} \pi_{k,l} \right)
                + \left(\sum_{k,l}\sum_{(i,j)\succ(k,l)} \pi_{i,j} \delta_{i,j}\right) \\
  \label{eq:obj_f_nested}
         & = \delta_{1,1}
                + \left(\sum_{k,l} \left[\delta_{k,l} \pi_{k,l} \ + \sum_{(i,j)\succ(k,l)} \pi_{i,j} \delta_{i,j}\right]\right)
\end{align}
which we refer to as the ``independent stages sum'' and ``nested stage sum'',
respectively \eqref{eq:obj_f_by_stage} and~\eqref{eq:obj_f_nested}.

%\begin{align} \label{eq:obj_f_decomposed}
%  f(\pi) & = \sum_{i,j} \pi_{i,j} d_{i,j} = \sum_{k,l} \sum_{(i,j)\succ(k,l)} \pi_{i,j} d_{i,j} \\
%         & = \sum_{k,l} \sum_{(i,j)\succ(k,l)} \pi_{i,j} \big(\delta_{1,1} + \delta_{k,l} + \delta_{i,j}\big) \\
%         & = \left(\sum_{k,l}\sum_{(i,j)\succ(k,l)} \pi_{i,j} \delta_{1,1}\right)
%                + \left(\sum_{k,l}\sum_{(i,j)\succ(k,l)} \pi_{i,j} \delta_{k,l}\right)
%                + \left(\sum_{k,l}\sum_{(i,j)\succ(k,l)} \pi_{i,j} \delta_{i,j}\right) \\
%         & = \delta_{1,1}
%                + \left(\sum_{k,l} \delta_{k,l} \sum_{(i,j)\succ(k,l)} \pi_{i,j} \right)
%                + \left(\sum_{k,l}\sum_{(i,j)\succ(k,l)} \pi_{i,j} \delta_{i,j}\right) \\
%         & = \delta_{1,1}
%                + \left(\sum_{k,l} \delta_{k,l} \pi_{k,l} \right)
%                + \left(\sum_{k,l}\sum_{(i,j)\succ(k,l)} \pi_{i,j} \delta_{i,j}\right) \\
%         & = \delta_{1,1}
%                + \left(\sum_{k,l} \left[\delta_{k,l} \pi_{k,l} \ + \sum_{(i,j)\succ(k,l)} \pi_{i,j} \delta_{i,j}\right]\right)
%\end{align}

\medskip

This motivates an alternative point of view on the optimization problem defining the Nested distance.
We introduce of a set of new variables $\pi_{k,l}$ and their corresponding constraints:
\[
  1 = \sum_{k,l} \pi_{k,l} \qquad \text{and} \qquad \forall k,l : \pi_{k,l} = \sum_{(i,j)\succ(k,l)} \pi_{i,j}.
\]
Of course, these ``new variables'' do not add any freedom to the feasible set
given the constraints they are subject to; this is only an affine lift of the original set.
Because of this, and of our notation, we still denote by $\pi$
the vector of variables with respect to which we're optimizing.

We can then see $\pi_{i,j}$ as decision variables for the third stage,
and $\pi_{k,l}$ for the second
\ (and, pushing the analogy further, $\pi_{1,1} = 1$ for the first one).

\subsection{Reformulating the constraints}
\label{sec:3stage_constr}

Now, to implement a recursive formulation of the problem,
we must also rewrite the constraints so that we can separate them more clearly.
This will be a particular application of the procedure from Pflug/Pichler
% TODO: cite (com seção e teo?)
that rewrites the Nested distance in dynamic programming form.

Let's focus on the constraints relative to $P$; the ones for $Q$ are analogous.
There are two groups of restrictions: the ones corresponding to the root node pair $(1,1)$,
and the ones corresponding to internal node pairs $(k,l)$.
Writing the sums explicitly, and using $\pi(i,j \mid k,l) = \frac{\pi_{i,j}}{\pi_{k,l}}$,
they are, respectively:
\begin{equation} \label{eq:constr_3stg}
  \begin{array}{ll}
    \sum_{j} \pi_{i, j} = P(i) & \text{for all $i$}\\
    \sum_{j \succ l} \pi_{i,j} = P(i \mid k) \pi_{k,l} & \text{for all $i \succ k$}
  \end{array}
\end{equation}
The second group of equations links the second to the third stage of the tree;
the first one, depending on the entire $\mathcal{B}$ tree, links the first to the third.
This is the one we must rewrite in order to obtain a relation only between the first and the second stage.
So, expand $\sum_j \pi_{i,j} = \sum_l \sum_{j \succ l} \pi_{i,j}$
and use the second group of equations (where $k$ must be the predecessor of $i$) to obtain
$P(i) = \sum_j \pi_{i,j} = \sum_l P(i \mid k) \pi_{k,l}$.
If we divide by $P(i \mid k)$ on both sides, and replace its definition, we get
\[
  P(k) = \sum_l \pi_{k,l}
\]
which does not include variables from the third stage anymore (no $i$ nor $j$),
and can be interpreted as only depending on the tree structure between the first and the second stages.

This shows that any feasible $\pi$ will also satisfy the simply coupled set of equations:
\begin{equation} \label{eq:constr_3stg_nested}
  \begin{array}{ll}
    \sum_l \pi_{k,l} = P(k) & \text{for all $k$} \\
    \sum_{j \succ l} \pi_{i,j} = P(i \mid k) \pi_{k,l} & \text{for all $i \succ k$}
  \end{array}
\end{equation}

The reverse inclusion can be obtained reversing the operations:
for each $i \succ k$, multiply the first constraint by (the constant) $P(i \mid k)$ on both sides,
and apply the second equation to arrive at $\sum_j \pi_{i,j} = P(k) P(i \mid k) = P(i)$.

\subsection{Recursive and independent formulations}

Now, we can replace our original 3-stage Nested Distance problem with the equivalent:
\begin{equation} \label{eq:ND_3stg_by_stage}
  \begin{array}{rlll}
    \displaystyle \min_{\pi}
    & \displaystyle \delta_{1,1} \ + {}
    & \displaystyle \sum_{k,l} \pi_{k,l} \delta_{k,l} \quad {} + {}
    & \displaystyle \sum_{k,l}\sum_{(i,j)\succ(k,l)} \pi_{i,j} \delta_{i,j} \\[1ex]
    \textrm{s.t.}
    && \text{for all nodes $k,l$, and} & \text{for all leaves $i \succ k$, $j \succ l$:} \\
    && \sum_{l} \pi_{k,l} = P(k) & \sum_{j \succ l} \pi_{i,j} = P(i \mid k) \pi_{k,l} \\
    && \sum_{k} \pi_{k,l} = Q(l) & \sum_{i \succ k} \pi_{i,j} = Q(j \mid l) \pi_{k,l} \\
    && \sum_{k,l} \pi_{k,l} = \pi_{1,1} = 1 & \sum_{(i,j)\succ(k,l)} \pi_{i,j} = \pi_{k,l}, \\
    && \pi_{k,l} \geq 0 & \pi_{i,j} \geq 0
  \end{array}
\end{equation}
There, we have separated in columns the components of the objective function,
and restrictions relative to the variables appearing on each column.
This shows that the symmetry of the structure of SWI tree
results in very similar sets of constraints.

By pushing the minimization inside the sum, we obtain both the recursive
and the independent formulations of the Nested distance.
The recursive one is independent of the tree structure, and splits the optimization problem
into a sum of optimization problems of choosing $\pi_{i,j}$ for $(i,j) \succ (k,l)$,
linked by the restrictions on $\pi_{k,l}$ from the second stage.
Concretely:
\begin{equation} \label{eq:ND_3stg_nested}
  \begin{array}{rll}
    \displaystyle \min_{\pi_{k,l}}
    & \displaystyle \delta_{1,1} \ + {}
    & \displaystyle \sum_{k,l} \pi_{k,l} \delta_{k,l} + \sum_{k,l} \Phi(k,l,\pi_{k,l}) \\[1ex]
    \textrm{s.t.}
    && \text{for all nodes $k,l$} \\
    && \sum_{l} \pi_{k,l} = P(k) \\
    && \sum_{k} \pi_{k,l} = Q(l) \\
    && \sum_{k,l} \pi_{k,l} = \pi_{1,1} = 1 \\
    && \pi_{k,l} \geq 0
  \end{array}
\end{equation}
where $\Phi(k,l,\pi_{k,l})$ is the optimal value of the optimization sub-problem:
\begin{equation} \label{eq:3stg_subprob}
  \begin{array}{rll}
    \displaystyle \min_{\pi_{i,j}} & \sum_{(i,j)\succ(k,l)} \pi_{i,j} \delta_{i,j} \\
    \textrm{s.t.}
      & \text{for all leaves $i \succ k$, $j \succ l$:} \\
      & \sum_{j \succ l} \pi_{i,j} = P(i \mid k) \pi_{k,l} \\
      & \sum_{i \succ k} \pi_{i,j} = Q(j \mid l) \pi_{k,l} \\
      & \sum_{(i,j)\succ(k,l)} \pi_{i,j} = \pi_{k,l} \\
      & \pi_{i,j} \geq 0
  \end{array}
\end{equation}
Observe that both the constraints and the objective function above are homogeneous of degree 1 on $\pi_{i,j}$
and $\pi_{k,l}$, so that $\Phi(k,l,\pi_{k,l}) = \pi_{k,l}\Phi(k,l,1)$.

Now, if the tree probabilities are stage-wise independent,
the subproblems for $\Phi(k,l,1)$ are all equal.
Indeed, given a pair of nodes $(k,l)$,
the leaves $(i,j)$ can be written as sequences $(k\otimes a, l\otimes b)$.
If $(k',l')$ is another pair of nodes from the same stage as $(k,l)$,
we can define $i' = k' \otimes a$ and $j' = l' \otimes b$,
so that the conditional probabilities $P(i \mid k)$ and $Q(j \mid l)$
are equal to $P(i' \mid k')$ and $Q(j' \mid l')$.
In the same way, $\delta_{i,j} = \delta_{i',j'}$,
since they are both equal to $d(a,b)$, the realization of the last stage.
This shows that all constants appearing in $\Phi(k,l,1)$
don't depend on $(k,l)$, so $\Phi(k,l,1) = \Phi$,
and we simplify further $\Phi(k,l,\pi_{k,l}) = \pi_{k,l} \Phi$.

Therefore, the objective function in~\eqref{eq:ND_3stg_nested} becomes
$\delta_{1,1} + \sum_{k,l} \pi_{k,l} \big(\delta_{k,l} + \Phi\big)$.
Since both $\delta_{1,1}$ and $\Phi$ don't depend
on the ``second-stage'' decision variables $\pi_{k,l}$,
the problem is equivalent to minimizing
$\sum_{k,l} \pi_{k,l} \delta_{k,l}$.

\medskip

This shows that the Nested Distance problem for stage-wise independent 3-stage trees
splits into two independent problems:
one for calculating the $\pi_{k,l}$, based on the first-to-second stage structure,
and one for calculating $\Phi$ based on the second-to-third-stage structure.
Each one corresponds to a Wasserstein distance calculation:
the first for the second-stage probabilities,
and the second for the third-stage probabilities.
Once both are solved, the actual probabilities $\pi_{i,j}$ can be obtained
by multiplying the corresponding optimal solution from $\Phi$
with the factor $\pi_{k,l}$ normalizing it for the $(k,l)$-subtree.

% vim:set spelllang=en:

%% file: multistage.tex
\section{Multi-stage setting}
\label{sec:multistage}

In this section, we generalize the result from the simple 3-stage setting
to compute the Nested Distance between two arbitrary SWI processes.

Our construction is based upon the successive equivalence of the original LP~\eqref{eq:ND}
with three other LP's, which transform the optimization problem in three different aspects.
The first one deals with the constraints,
and the second introduces the Benders decomposition / dynamic programming.
% TODO: following Pflug/Pichler (maybe Pichler Kovacevich?)
An interesting aspect of this decomposition is that each subproblem in the dynamic programming formulation
is also a Nested Distance between certain subtrees.
From this formulation, the main result of this article is straightforward,
coming as the third LP, where one takes advantage of the independence of the stages to reach further simplification.

% In fact, we can use 1-Markov or even k-Markov and obtain (some) simplification

\subsection{Rewriting constraints in successor form}

Similarly to the 3-stage case, we introduce the notation
$\pi_{k,l}$ for the subtree probability $\sum_{(i,j) \succ (k,l)} \pi_{i,j}$
ranging over all descendants from $k$ and $l$.
Replacing the conditional probability $\pi(i,j|k,l)$ by the ratio of joint probabilities $\pi_{i,j}/\pi_{k,l}$,
the Nested Distance LP, \eqref{eq:ND}, becomes
\begin{equation} \label{eq:ND_LP-LV}
\begin{array}{rll}
\displaystyle \min_\pi  & \displaystyle \sum_{i,j}\pi_{i,j} \cdot d_{i,j}^p \\
\textrm{s.t.} & \text{for all nodes $k \in \mathcal{A}_t$ and $l \in \mathcal{B}_t$:} \\
              & \sum_{j \succ l} \pi_{i, j} = P(i \mid k) \pi_{k,l}
              & \text{for all $i \succ k$}, \\
              & \sum_{i \succ k} \pi_{i, j} = Q(j \mid l) \pi_{k,l}
              & \text{for all $j \succ l$}, \\[1ex]
              & \sum_{i,j} \pi_{i,j} = 1, \\
              & \pi_{i,j} \geq 0,
              % & \text{for all leaves $i \in \mathcal{A}_T$ and $j \in \mathcal{B}_T$}.
\end{array}
\end{equation}
In what follows, we will omit the statement about the nodes $k$ and $l$
being nodes belonging to stage $t$.
Unless stated otherwise, each constraint involving nodes $k$ and $l$ is imposed
for any pair $(k,l)$ of nodes from $\mathcal{A}_t\times\mathcal{B}_t$,
for all stages $t$ ranging from $1$ up to $T$.

\medskip

If we now consider $\pi_{k,l}$ as variables, and introduce the constraints corresponding to their definition,
we obtain a problem equivalent to \eqref{eq:ND_LP-LV} that helps with the transition to the dynamic programming form:
\begin{equation*}
\begin{array}{rll}
\displaystyle \min_\pi  & \displaystyle \sum_{i,j}\pi_{i,j} \cdot d_{i,j}^p \\
\textrm{s.t.} & \sum_{j \succ l} \pi_{i, j} = P(i \mid k) \pi_{k,l}
              & \text{for all $i \succ k$}, \\
              & \sum_{i \succ k} \pi_{i, j} = Q(j \mid l) \pi_{k,l}
              & \text{for all $j \succ l$}, \\[1ex]
              & \sum_{(i,j) \succ (k,l)} \pi_{i,j} = \pi_{k,l}, \\[1ex]
              & \sum_{i,j} \pi_{i,j} = 1, \\
              & \pi_{i,j} \geq 0, \pi_{k,l} \geq 0.
\end{array}
\end{equation*}

Better still, for this purpose, would be rewriting the conditional marginal distribution and sum-to-one constraints
in terms of of successors nodes rather than leaf nodes, analogously to section~\ref{sec:3stage_constr}:
\begin{equation}
\begin{array}{rll} \label{eq:ND_LP-SC}
\displaystyle \min_\pi  & \displaystyle \sum_{i,j}\pi_{i,j} \cdot d_{i,j}^p \\
\textrm{s.t.} & \sum_{s \in l_+} \pi_{r, s} = P(r \mid k) \pi_{k,l}
& \text{for all $r \in k_+$}, \\
& \sum_{r \in k_+} \pi_{r, s} = Q(s \mid l) \pi_{k,l}
& \text{for all $s \in l_+$}, \\[1ex]
              & \sum_{(r,s) \in (k,l)_+} \pi_{r,s} = \pi_{k,l}, \\
              & \pi_{1,1} = 1, \\
& \pi_{k,l} \geq 0.
\end{array}
\end{equation}
This is a valid transformation, as shown by the following lemma:
\begin{lemma}
  The following two sets of constraints are equivalent:
  \begin{equation}\label{eq:constr_deep}
    \sum_{j \succ l} \pi_{i, j} = P(i \mid k) \pi_{k,l} \quad\text{for all $k,l$}
  \end{equation}
  \begin{equation}\label{eq:constr_sw}
    \sum_{s \in l+} \pi_{r,s} = P(r \mid k) \pi_{k,l} \quad\text{for all $k,l$, $r \in k+$}
  \end{equation}
\end{lemma}
\begin{proof}

For all $r \in k_+$ and $i \succ r$ we have, starting from~\eqref{eq:constr_deep}, reversed:
\begin{equation}
P(i \mid k) \pi_{k,l} = \sum_{j \succ l} \pi_{i, j} = \sum_{s \in l_+} \sum_{j \succ s} \pi_{i, j} = \sum_{s \in l_+} P(i \mid r) \pi_{r,s}.
\end{equation}
where we applied again~\eqref{eq:constr_deep}, but with $(r,s)$ instead of $(k,l)$ on the last step.
Therefore,
\begin{equation}
  \sum_{s \in l_+} \pi_{r,s}
  = \frac{P(i \mid k)}{P(i \mid r)} \pi_{k,l}
  = P(r \mid k) \pi_{k,l},
\end{equation}
for all $k \in \mathcal{A}_t$, $l \in \mathcal{B}_t$ and all $r \in k_+$.
This shows that~\eqref{eq:constr_deep} implies~\eqref{eq:constr_sw}.

For the other direction, for all $k \in \mathcal{A}_t$, $l \in \mathcal{B}_t$ and all $i \succ k$ we have,
applying~\eqref{eq:constr_sw} at all nodes $r,m,\ldots$ in the path from $k$ to $i$,
and corresponding nodes $s,n,\ldots$ from $l$:
\begin{align*}
\pi_{k,l} & = \frac{P(k)}{P(r)} \sum_{s \in l_+} \pi_{r,s}
            = \frac{P(k)}{P(r)} \frac{P(r)}{P(m)} \sum_{s \in l_+} \sum_{n \in s_+} \pi_{m,n} \\
          & = \cdots = \frac{P(k)}{P(i)} \sum_{j \succ l} \pi_{i,j}.
\end{align*}
Therefore, those constraints are equivalent.
\end{proof}

%\begin{equation*}
%\begin{array}{rl}
%\displaystyle \min_\Pi  & \displaystyle \sum_{i,j}\pi_{i,j} \cdot d_{i,j}^p \\
%\textrm{s.t.} & (RP) + (RQ) + (Pos) + (Norm)
%\end{array}
%\end{equation*}

\subsection{Dynamic programming}

We note that \eqref{eq:ND_LP-SC} has a recursive constraint structure regarding the probability $\pi_{k,l}$ of a given node pair and the probabilities $\pi_{r,s}$ of its successor node pairs.
In order to define our dynamic programming problem, we should state the objective function in a similar form.
A recursive representation of the objective function is obtained by decomposing the sum over all leaves in a nested sum over successor nodes from root to the leaves:
\begin{equation} \label{eq:obj_LP-SC}
\sum_{i,j} \pi_{i,j} \cdot d_{i,j}^p = \sum_{(k_2,l_2) \in (1_\mathcal{A},1_\mathcal{B})_+} \cdots \sum_{(i,j) \in (k_{T-1},l_{T-1})_+} \pi_{i,j} \cdot d_{i,j}^p.
\end{equation}
From representation \eqref{eq:obj_LP-SC}, we can define recursively a function $\Phi_t(k,l,\pi_{k,l})$ that describes our first dynamic programming formulation of the Nested Distance.
For the last stage $T$, we define $\Phi_{T}(i,j,\pi_{i,j})$ as $\pi_{i,j} \cdot d_{i,j}^p$ and, for a general stage $t$,
we define $\Phi_t(k,l,\pi_{k,l})$ as the optimal value of the following optimization problem
\begin{equation} \label{eq:phi_klm}
\begin{array}{rll}
\displaystyle \min_\pi  & \displaystyle \sum_{(r,s) \in (k,l)_+}  \Phi_{t+1}(r,s,\pi_{r,s}) \\
\textrm{s.t.} & \sum_{s \in l_+} \pi_{r, s} = P(r \mid k) \pi_{k,l}
& \text{for all $r \in k_+$}, \\
& \sum_{r \in k_+} \pi_{r, s} = Q(s \mid l) \pi_{k,l}
& \text{for all $s \in l_+$}, \\[1ex]
              & \sum_{(r,s) \in (k,l)_+} \pi_{r,s} = \pi_{k,l}, \\
& \pi_{r,s} \geq 0,
\end{array}
\end{equation}
where both the node pair $(k,l) \in \mathcal{A}_t \times \mathcal{B}_t$
and the total weight $\pi_{k,l}$ are \emph{fixed}.

Note that $\Phi_t(k,l,\pi_{k,l})$ is positive homogeneous in $\pi_{k,l}$, that is,
\begin{equation} \label{eq:pos-hom}
\Phi_t(k,l,\pi_{k,l}) = \pi_{k,l} \cdot \Phi_t(k,l,1),
\end{equation}
for any non-negative number $\pi_{k,l}$.
This may be proved by (backwards) induction on the stage $t$.
By definition, \eqref{eq:pos-hom} is true for the last stage $T$.
Now, since the constraints are homogeneous on $(\pi_{k,l}, \pi_{r,s})$,
the feasible set for $\Phi_t(k,l,\alpha\pi_{k,l})$ is just a scaling of the one for $\Phi_t(k,l,\pi_{k,l})$.
From the induction hypothesis, the objective function is positive homogeneous:
\begin{equation} \label{eq:pos-hom_sum}
\sum_{(r,s) \in (k,l)_+}  \Phi_{t+1}(r,s,\pi_{r,s}) = \sum_{(r,s) \in (k,l)_+} \pi_{r,s} \Phi_{t+1}(r,s,1),
\end{equation}
so we also get a scaling between the \emph{solutions} in each case,
and their optimal value.
This shows that the sub-problem for stage $t$ is positive homogeneous.

\medskip

We assert that $\Phi_{1}(\mathbf{1}_{\mathcal{A}},\mathbf{1}_{\mathcal{B}},1)$
is the Nested Distance between probability trees $\mathcal{A}$ and $\mathcal{B}$.
Indeed, if we replace recursively the definition of $\Phi_t(k,l,\pi_{k,l})$,
group the nested minimization in a single minimization problem
and use the identity \eqref{eq:obj_LP-SC} for the objective function,
then we get back formulation \eqref{eq:ND_LP-SC}.

% Explicar melhor, 3 estagios
\bigskip

Finally, we can give a second dynamic programming formulation of the Nested distance.
Define $\mathbf{d}^p(k,l) := \Phi_t(k,l,1)$, a function we refer to as the \emph{sub-tree distance} between $k$ and $l$.
It is instructive to note that $\mathbf{d}^p(i,j)$ is equal to $d_{i,j}^p$, by definition of $\Phi_T(i,j,1)$, and $\mathbf{d}^p(k,l)$ is equal to the optimization problem
\begin{equation} \label{eq:ND_BD}
\begin{array}{rll}
\displaystyle \min_\pi  & \displaystyle \sum_{(r,s) \in (k,l)_+} \pi_{r,s} \cdot \mathbf{d}^p(r,s) \\
\textrm{s.t.} & \sum_{s \in l_+} \pi_{r, s} = P(r \mid k)
& \text{for all $r \in k_+$}, \\
& \sum_{r \in k_+} \pi_{r, s} = Q(s \mid l)
& \text{for all $s \in l_+$}, \\[1ex]
              & \sum_{(r,s) \in (k,l)_+} \pi_{r,s} = 1, \\
& \pi_{r,s} \geq 0,
\end{array}
\end{equation}
by definition \eqref{eq:phi_klm} of $\Phi_t(k,l,1)$,
the representation \eqref{eq:pos-hom_sum} of the objective function
and backward induction on the stage that $(k,l)$ belongs to.
Therefore, $\mathbf{d}(\mathbf{1}_\mathcal{A},\mathbf{1}_\mathcal{B})$ is the Nested Distance between probability trees $\mathcal{A}$ and $\mathcal{B}$.
We emphasize that both dynamic programming formulation \eqref{eq:phi_klm} and \eqref{eq:ND_BD} are equivalent due to the positive homogeneity property of $\Phi_t(k,l,\pi_{k,l})$.

\medskip

An interpretation of $\mathbf{d}(k,l)$ can be given in terms of the Nested Distance between certain sub-trees,
which explains its name.
Let $\mathcal{A}(k)$ be the probability tree defined as the combination between a path from the root node $\mathbf{1}_\mathcal{A}$ to the node $k$ and the subprobability tree rooted at $k$.
Let $\mathcal{B}(l)$ be defined analogously.
We claim that $\mathbf{d}(k,l)$ is equal to the Nested Distance between $\mathcal{A}(k)$ and $\mathcal{B}(l)$, i.e.,
\[
  \mathbf{d}(k,l) = \mathbf{d}_{\textrm{N}}(\mathcal{A}(k),\mathcal{B}(l)).
\]
Indeed, for stages at or after $t$, the optimization problems for $\mathcal{A}(k)\times\mathcal{B}(l)$ and $\mathcal{A}\times\mathcal{B}$
are the same,
and for stages prior to $t$,
the optimization problem for $\mathbf{d}_{\textrm{N}}$ is trivial since there's only one successor node pair.

\subsection{Stagewise independence}

We are now in position to prove the main result of this part of the paper.
Notice that, while the dynamic programming form expressed in equation~\eqref{eq:ND_BD}
is valid for arbitrary trees,
the final simplification we get here comes from two different aspects of stagewise independence:
the probability structure on the stochastic process and the decomposition of the distance function.

\begin{thm}
Let $\mathcal{A}$ and $\mathcal{B}$ be two stagewise independent trees
and $d(\xi,\zeta)$ be a distance function given by the weighted sum of distances between coordinates $\xi_t$~and~$\zeta_t$:
\[
d^p(\xi,\zeta) := w_1 \cdot d_1^p(\xi_1,\zeta_1) + \dots + w_T \cdot d_T^p(\xi_T,\zeta_T).
\]
Then the Nested Distance between $\mathcal{A}$ and $\mathcal{B}$ is equal to the weighted sum of the Wassertein distances between the marginal distribution of each stage:
\begin{equation} \label{eq:ND-WS}
\mathbf{d}_\textrm{N}^p(\mathcal{A},\mathcal{B}) = w_{1}\cdot\mathbf{d}_\mathrm{W}^p(P_{1}, Q_{1}) + \dots + w_{T}\cdot\mathbf{d}_\mathrm{W}^p(P_{T}, Q_{T}).
\end{equation}
\end{thm}

\begin{proof}
Under the same hypothesis, we will prove
% by induction from last to first stage
a more general result about the sub-tree distance $\mathbf{d}^p(k,l)$
which includes \eqref{eq:ND-WS} as a particular case.

Let $k$ and $l$ be two nodes from stage $t$ belonging to $\mathcal{A}$ and $\mathcal{B}$, respectively.
Say $k$ corresponds to all processes that start with $(\xi_1^k, \xi_2^k, \ldots, \xi_t^k)$,
and $l$ to processes starting with $(\zeta_1^l, \zeta_2^l, \ldots, \zeta_t^l)$.
We claim that the sub-tree distance between $k$ and $l$ is equal to
the weighted sum of distances between the ``past'' outcomes $\xi_\tau^k$ and $\zeta_\tau^l$
up to stage $t$,
plus the Wassertein distances between the marginal distributions
from stage $t+1$ up to the last stage $T$:
\begin{align}
\begin{split} \label{eq:ND-ND}
\mathbf{d}^p(k,l) = &  w_1 \cdot d_1^p\left(\xi_1^k, \zeta_1^l\right) + \dots + w_t \cdot d_t^p\left(\xi_t^k, \zeta_t^l\right) \\
& + w_{t+1}\cdot\mathbf{d}_\mathrm{W}^p(P_{t+1}, Q_{t+1}) + \dots + w_{T}\cdot\mathbf{d}_\mathrm{W}^p(P_{T}, Q_{T}).
\end{split}
\end{align}
The statement \eqref{eq:ND-ND} is equivalent to \eqref{eq:ND-WS} when the nodes $k$~and~$l$ are the root nodes $\mathbf{1}_{\mathcal{A}}$~and~$\mathbf{1}_{\mathcal{B}}$, respectively.

\newcommand{\R}{\mathrm{R}}
In order to simplify notation, we denote by $\R_{t+1}$ the sum of Wasserstein distances
%between marginal distribution
from stage $t+1$ up to the last stage $T$:
\[
\R_{t+1} := w_{t+1}\cdot\mathbf{d}_\mathrm{W}^p(P_{t+1}, Q_{t+1}) + \dots + w_{T}\cdot\mathbf{d}_\mathrm{W}^p(P_{T}, Q_{T}).
\]
We proceed by backward induction over the stage $t$.
Indeed, the identity \eqref{eq:ND-ND} is trivial for the last stage $T$, since it is the definition of $\mathbf{d}^p(i,j)$:
\[
  \mathbf{d}^p(i,j) = d^p_{i,j} = d^p(\xi^i, \zeta^j).
\]
By the induction hypothesis, equation~\eqref{eq:ND-ND} holds for all successor pairs $(r,s)$
of any node pair $(k,l)$ from stage $t$:
\begin{equation}\label{eq:ND-ND_succ}
\begin{aligned}
\mathbf{d}^p(r,s) & = w_1 \cdot d_1^p\left(\xi_1^r, \zeta_1^s\right) + \dots + w_{t+1} \cdot d_{t+1}^p\left(\xi_{t+1}^r, \zeta_{t+1}^s\right) \\
& \qquad + \R_{t+2}.
%& + w_{t+2}\cdot\mathbf{d}_\mathrm{W}^p(P_{t+2}, Q_{t+2}) + \dots + w_{T}\cdot\mathbf{d}_\mathrm{W}^p(P_{T}, Q_{T}).
\end{aligned}
\end{equation}
Now, since the outcomes $\xi_\tau^r$ and $\zeta_\tau^s$ for $\tau \leq t$
are exactly the outcomes $\xi_\tau^k$ and $\zeta_\tau^l$ given by $k$ and $l$,
we can rewrite equation~\eqref{eq:ND-ND_succ} to express the dependence on each node explicitly:
\begin{equation} \label{eq:ND-ND_succ2}
\begin{aligned}
\mathbf{d}^p(r,s) & =  w_1 \cdot d_1^p\left(\xi_1^k, \zeta_1^l\right) + \dots + w_{t} \cdot d_{t}^p\left(\xi_{t}^k, \zeta_{t}^l\right)  \\
& \quad + w_{t+1} \cdot d_{t+1}^p\left(\xi_{t+1}^r, \zeta_{t+1}^s\right) + \R_{t+2}.
\end{aligned}
\end{equation}
Note that the only element of \eqref{eq:ND-ND_succ2}  depending on the successor node $(r,s)$ is the distance $d_{t+1}^p\left(\xi_{t+1}^r, \zeta_{t+1}^s\right)$ between a fixed realization from stage $t+1$.
By definition~\eqref{eq:ND_BD} of sub-tree distance and from equation \eqref{eq:ND-ND_succ2}, we have that $\mathbf{d}^p(k,l)$ is equal to the optimal value of the optimization problem
\begin{equation} \label{eq:ND_BD-stgind}
\begin{array}{rll}
\displaystyle \min_\pi  & \displaystyle \sum_{(r,s) \in (k,l)_+} \pi_{r,s} \cdot  w_{t+1}d_{t+1}^p\left(\xi_{t+1}^r, \zeta_{t+1}^s\right) \\
\textrm{s.t.} & \sum_{s \in l_+} \pi_{r, s} = P(r \mid k)
& \text{for all $r \in k_+$}, \\
& \sum_{r \in k_+} \pi_{r, s} = Q(s \mid l)
& \text{for all $s \in l_+$}, \\[1ex]
& \sum_{(r,s) \in (k,l)_+} \pi_{r,s} = 1, \\
& \pi_{r,s} \geq 0,
\end{array}
\end{equation}
plus $w_1 \cdot d_1^p\left(\xi_1^k, \zeta_1^l\right) + \dots + w_{t} \cdot d_{t}^p\left(\xi_{t}^k, \zeta_{t}^l\right) + \R_{t+2}$.
But the optimal value of \eqref{eq:ND_BD-stgind} is the Wassertein distance $\mathbf{d}_{\mathrm{W}}(P_t,Q_t)$ between marginal distributions $P_t$ and $Q_t$,
which concludes the induction.
\end{proof}

% vim:set spelllang=en:

%% file: example.tex
\section{Example}

We can use the decomposed formula of the Nested distance in the stagewise independent setting
to simplify the problem of optimal scenario generation.
This is suitable for practical applications where the number of stages and,
consequently, scenarios are huge.
Indeed, it is very common that both the original process $P$
and the approximant $Q$ must be stagewise independent,
for instance in order to employ the SDDP algorithm.

By the SWI hypothesis, the scenario probability
$P\big([\xi_1,\xi_2, \ldots, \xi_T]\big)$
decomposes as a product of stagewise probabilities, $\prod P_i(\xi_i)$,
and that also follows for~$Q$.
This implies that, for all such $Q$, the Nested distance $d_{ND}(P,Q)$ is the
sum of Wasserstein distances between all marginal distributions pairs,
$\sum_{t} d_W(P_t,Q_t)$.
This splits the optimization problem of optimal scenario generation
into $T$ independent subproblems,
\[
\min_{Q \in \mathrm{SWI}} d_{ND} (P,Q) = \sum_{t = 1}^{T} \min_{Q_i} d_{W} (P_i,Q_i)
\]
where each $d_{W} (P_i,Q_i)$ is faster to compute and
the whole procedure is easily parallelizable.

In the particular case of a SWI process $P$ with equally probable scenarios and
a Nested distance induced by a stagewise quadratic (euclidean) distance
the minimizer of the Wasserstein distance can be calculated via the K-means
algorithm~\cite[section 3]{KovacevicP15}.
It is worth mentioning that such approach is employed nowadays in the Brazilian Mid
and Long Term Power System models \cite{penna2009definiccao}, however, as far
as we know there is no formalization of it.